%% file: nonlinear_fractional_diffusion-2.tex
\documentclass{amsproc}

\usepackage{amssymb}

\usepackage{}
\usepackage{mathrsfs,mathtools}
\usepackage{enumerate}
\usepackage{enumitem}
\usepackage{cleveref}
\usepackage{autonum}
\usepackage{array,multirow}

\newtheorem{theorem}{Theorem}[section]

\newtheorem{proposition}[theorem]{Proposition}

\theoremstyle{definition}

\newtheorem{example}[theorem]{Example}

\theoremstyle{remark}
\newtheorem{remark}[theorem]{Remark}

\numberwithin{equation}{section}

\def\l{\left}
\def\r{\right}

\newcommand{\SZ}{\Pi_h}

\def\mR{\mathbb{R}}
\def\mRd{{\mathbb{R}^d}}

\def\PV{\text{P.V.}}
\def\Gts{\widetilde{G}_s}
\def\Th{\mathcal{T}_h}

\def\Cl{\Pi_h}
\def\ve{\varepsilon}
\newcommand{\eps}{\varepsilon}
\def\pO{\partial\Omega}

\newcommand{\phii}{\varphi}
\newcommand{\pp}{\partial}

\def\dist{\textrm{dist}}
\def\oT{T}
\def\oTp{T'}

\DeclareMathOperator{\supp}{supp}

\newcommand{\Laps}{(-\Delta)^s}

\newcommand{\F}{\mathcal{F}_s}

\newcommand{\Fou}{\mathscr{F}}
\newcommand{\argmin}{\textrm{argmin}}

\newcommand{\tH}{{\widetilde H}}
\newcommand{\Hs}{H^s(\Omega)}
\newcommand{\tHs}{{\widetilde H}^s(\Omega)}
\newcommand{\R}{{\mathbb{R}}}

\newcommand{\T}{{\mathcal{T}_h}}
\newcommand{\K}{{\mathcal{K}}}

\newcommand{\gain}{\gamma}
\newcommand{\gainB}{t}

\newcommand{\vA}{\mathbf{A}}
\newcommand{\vU}{\mathbf{U}}
\newcommand{\vF}{\mathbf{F}}
\newcommand{\vB}{\mathbf{B}}

\newcommand{\x}{\texttt{x}}

\usepackage[usenames,dvipsnames]{color}

\newcommand{\rhn}[1]{{\color{black}#1}}

\begin{document}

\title[Linear and Nonlinear Fractional Elliptic Problems]{\rhn{Linear and Nonlinear Fractional Elliptic Problems}}


\author[J.P.~Borthagaray]{Juan Pablo~Borthagaray}
\address[J.P.~Borthagaray]{Department of Mathematics, University of Maryland, College Park, MD 20742, USA}
\email{jpb@math.umd.edu}
\thanks{JPB has been supported in part by NSF grant DMS-1411808}

\author[W.~Li]{Wenbo~Li}
\address[W.~Li]{Department of Mathematics, College Park, MD 20742, USA}
\email{wenboli@math.umd.edu}
\thanks{WL has been supported in part by NSF grant DMS-1411808 and the Patrick and Marguerite Sung Fellowship in Mathematics.}

\author[R.H.~Nochetto]{Ricardo H.~Nochetto}
\address[R.H.~Nochetto]{Department of Mathematics and Institute for Physical Science and Technology, University of Maryland, College Park, MD 20742, USA}
\email{rhn@math.umd.edu}
\thanks{RHN has been supported in part by NSF grant DMS-1411808}

\date{}

\begin{abstract}
This paper surveys recent analytical and numerical research on linear problems for the integral fractional Laplacian, fractional obstacle problems, and fractional minimal graphs. The emphasis is on the interplay between regularity, including boundary behavior, and approximability by piecewise linear finite element methods. We discuss several error estimates on graded meshes, and computational challenges associated to implementing and solving efficiently the ensuing integral equations, along with numerical experiments.
\end{abstract}

\maketitle

\section{Introduction} \label{sec:intro}

Diffusion, which is one of the most common physical processes, is the net movement of particles from a region of higher concentration to a region of lower concentration. The assumption that particles respond to Brownian motion leads to classical models of diffusion \cite{Albert}, that have been well studied for a long time. 
Fick's first law states that the magnitude of the diffusive flux is proportional to the concentration gradient; by now, it is clear that such a constitutive relation is a questionable model for numerous phenomena \cite{MetzlerKlafter}. When the associated underlying stochastic process is not given by Brownian motion, the diffusion is regarded as \emph{anomalous}. In particular, anomalous superdiffusion refers to situations that can be modeled using fractional spatial derivatives or fractional spatial differential operators.

Integer-order differentiation operators are local because the derivative of a function at a given point depends only on the values of the function in an infinitesimal neighborhood of it. In contrast, fractional-order derivatives are nonlocal, integro-differential operators. A striking example of such an operator is the fractional Laplacian of order $s \in (0,1)$, which we will denote by $(-\Delta)^s$, and is given by
\begin{equation} \label{eq:def_Laps}
(-\Delta)^s u (x) := C_{d,s} \, \PV \int_\mRd \frac{u(x)-u(y)}{|x-y|^{d+2s}} \; dy , \quad C_{d,s} := \frac{2^{2s} s \Gamma(s+\frac{d}{2})}{\pi^{d/2} \Gamma(1-s)} .
\end{equation}
We refer to \cite{Valdinoci} for an illustration of how the heat equation involving the fractional Laplacian arises from a simple random walk with jumps. The nonlocal structure of the operator \eqref{eq:def_Laps} is apparent: to evaluate $(-\Delta)^s u$ at a spatial point, information involving all spatial points is needed.

This work deals with fractional diffusion. Our main goal is to review finite element methods (FEMs) to approximate solution of elliptic problems involving $(-\Delta)^s$ or related operators on bounded domains. We shall not discuss methods for the {\em spectral} fractional Laplacian; the surveys \cite{BBNOS18,Lischke_et_al18} offer comparison between such an operator and the fractional Laplacian \eqref{eq:def_Laps} and review other numerical methods. 
\rhn{We point out that the fractional Laplacian \eqref{eq:def_Laps} of order $s \in (0,1)$ is the infinitesimal generator of a $2s$-stable L\'evy process. In this regard, problems on a bounded domain with homogeneous Dirichlet boundary conditions arise when the process is killed upon exiting the domain.
}

Throughout this work, we assume that $\Omega$ is a bounded \rhn{Lipschitz domain. Whenever additional assumptions on $\pp\Omega$ are required, we shall state them explicitly. 
Even though there is a wide variety of numerical methods for fractional-order problems available in the literature \cite{Lischke_et_al18}, in this work we shall focus on piecewise linear finite element methods as in \cite{BBNOS18}. We emphasize the interplay between regularity, including boundary behavior, and approximability. In fact, the convergence rates achievable for the fractional elliptic PDEs discussed below, both linear and nonlinear,  are limited by the presence of an algebraic boundary layer regardless of the regularity of $\partial\Omega$ and the polynomial degree for shape regular elements.
}

The paper is organized as follows. \Cref{sec:linear} deals with the homogeneous Dirichlet problem for the fractional Laplacian in $\Omega$; we discuss regularity of solutions and discuss both theoretical and computational aspects of conforming finite element discretizations. We also comment on some recent applications of this approach and on an alternative 
nonconforming FEM based on a Dunford-Taylor representation of the weak form of $\Laps$. 
Afterwards, in \Cref{sec:obstacle} we address the obstacle problem for the fractional Laplacian. To derive optimal convergence estimates, we focus on weighted Sobolev regularity, where the weight is a power of the distance to the boundary of $\Omega$. These estimates follow from a precise quantification of boundary regularity of solutions and how solutions detach from the obstacle. Finally, \Cref{sec:NMS} deals with fractional minimal graphs, which in fact are subgraphs that minimize a suitable nonlocal perimeter. We formulate a variational form for this problem, which is nonlinear and degenerate. We report on approximation properties of a conforming finite element scheme, and show convergence rates with respect to a novel geometric quantity. The paper concludes with a couple of computational explorations of the behavior of fractional minimal graphs for $d=2$.


\input{linear}
%
\input{obstacle}
%
\input{minimal_surfaces}

\bibliography{NFD}
\bibliographystyle{amsplain}

\end{document}

%% file: linear.tex
\section{Linear Problems} \label{sec:linear}
In this section we consider the homogeneous Dirichlet problem for the fractional Laplacian \eqref{eq:def_Laps}. Given $f : \Omega \to \mR$,  one seeks a function $u$ such that
	\begin{equation}\label{eq:Dirichlet}
		\left\lbrace\begin{array}{rl}
			(-\Delta)^s u = f &\mbox{ in }\Omega, \\
				u = 0 &\mbox{ in }\Omega^c :=\mRd\setminus\Omega.
		\end{array} \right.
\end{equation}

\subsection{Variational Formulation}\label{S:variational_form}
%
The natural variational framework for \eqref{eq:Dirichlet} is within the fractional Sobolev space $\tHs$, that is defined by
\[
\tHs := \{ v \in H^s(\mRd) \colon \supp v \subset \overline{\Omega} \} .
\]
We refer to \cite{BoCi19} for definitions and elementary properties of fractional Sobolev spaces. Here we just state that on the space $\tHs$, because of the fractional Poincar\'e inequality, the natural inner product is equivalent to
\begin{equation} \label{eq:defofinnerprod}
  (v,w)_s := \frac{C_{d,s}}2 \iint_{Q_\Omega} \frac{(v(x)-v(y))(w(x) - w(y))}{|x-y|^{d+2s}} \; dx \; dy,
\end{equation}
where $Q_\Omega = (\mRd\times\mRd) \setminus(\Omega^c\times\Omega^c)$.
The corresponding norm is $\|v\|_{\tHs} := (v,v)_s^{1/2}$.

The duality pairing between $\tHs$ and its dual $H^{-s}(\Omega)$ shall be denoted by $\langle \cdot , \cdot \rangle$. In view of \eqref{eq:defofinnerprod} we see that, whenever $v \in \tHs$ then $(-\Delta)^s v \in H^{-s}(\Omega)$ and that
\begin{equation}
\label{eq:innerprodisduality}
  (v,w)_s = \langle (-\Delta)^s v, w \rangle, \quad \forall w \in \tHs.
\end{equation}
Therefore, given $f \in H^{-s}(\Omega)$, the weak formulation of \eqref{eq:Dirichlet} reads: find $u \in \tHs$ such that
\begin{equation} \label{eq:weak_linear}
(u, v)_s = \langle f, v \rangle \quad \forall v \in \tHs.
\end{equation}
Existence and uniqueness of weak solutions, and stability of the solution map $f \mapsto u$, follow straightforwardly from the Lax-Milgram Theorem.

\subsection{Regularity} \label{sec:regularity_linear}
A priori, it is not clear how smooth weak solutions are. If $f$ is more regular than $H^{-s}(\Omega)$, then $u$ could be expected to be more regular than $\tHs$. 
We now review some results regarding regularity of solutions to problem \eqref{eq:weak_linear}. Since our main interest is to derive convergence rates for finite element discretizations, we shall focus on Sobolev regularity estimates.

Recently, using Fourier analytical tools, Grubb \cite{Grubb} obtained estimates of solutions in terms of the so-called H\"ormander $\mu$-spaces \cite{Hormander}, but such estimates can be reinterpreted in terms of standard Sobolev spaces. A drawback of the following result from \cite{Grubb} is that it assumes the domain $\Omega$ to have smooth boundary, which is a too restrictive condition for finite element applications.

\begin{theorem}[regularity on smooth domains] \label{T:reg_grubb}
Let $\Omega$ be a domain with $\pp\Omega \in C^\infty$, $s \in (0,1)$,  $f \in H^r(\Omega)$ for some $r\ge -s$, $u$ be the solution of \eqref{eq:weak_linear} and $\gain = \min \{ s + r, 1/2 -\eps \}$, with $\eps > 0$ arbitrarily small. Then, $u \in \widetilde{H}^{s + \gain}(\Omega)$ and the following regularity estimate holds: 
\[
  \| u  \|_{\widetilde{H}^{s + \gain}(\Omega)} \le C(\Omega,d,s,\gain) \|f \|_{H^r(\Omega)}.
\]
\end{theorem}

A rather surprising feature of the previous result is that no matter how smooth the right hand side $f$ is, we cannot guarantee that solutions are any smoother than $\widetilde{H}^{s + 1/2 - \eps}(\Omega)$. Indeed, because the fractional Laplacian is an operator of order $2s$, it could be expected to have a lift of order $2s$. As the following example \cite{Getoor} shows, such a reasoning is incorrect, and \Cref{T:reg_grubb} is sharp.

\begin{example}[limited regularity]\label{ex:nonsmooth}
Consider $\Omega = B(0,1) \subset \mRd$ and $f \equiv 1$. Then, the solution to \eqref{eq:Dirichlet} is given by 
\begin{equation} \label{eq:getoor}
u(x) = \frac{\Gamma(\frac{d}{2})}{2^{2s} \Gamma(\frac{d+2s}{2})\Gamma(1+s)} ( 1- |x|^2)^s_+,
\end{equation}
where $t_+ =\max\{t,0\}$.
\end{example}

As \Cref{ex:nonsmooth} illustrates, rough boundary behavior causes the reduced Sobolev regularity of solutions. Ros-Oton and Serra \cite{RosOtonSerra} studied problem \eqref{eq:Dirichlet} using potential theory tools, and were able to obtain a fine characterization of boundary behavior of solutions, that led them to deduce H\"older regularity estimates. In particular, it turns out that the asymptotic expansion
\begin{equation} \label{eq:boundary_behavior}
u (x) \approx d(x, \pp\Omega)^s \varphi(x),
\end{equation}
where $\varphi$ is a smooth function, is generic.

The H\"older estimates from \cite{RosOtonSerra} give rise to Sobolev estimates for solutions in terms of H\"older norms of the data.  To capture the boundary behavior, reference \cite{AcosBort2017fractional} introduced fractional weighted norms, where the weight is a power of the distance to the boundary, and developeds estimates in such norms.
We denote
\[
\delta(x,y) := \min \big\{ \dist(x, \pp \Omega), \dist(y, \pp \Omega) \big\}
\]
and, for $\ell = k + s$, with $k \in \mathbb{N}$ and $s \in (0,1)$, and $\kappa \ge 0$, we define the norm
\[
\| v \|_{H^{\ell}_\kappa (\Omega)}^2 := \| v \|_{H^k (\Omega)}^2 + \sum_{|\beta| = k } 
\iint_{\Omega\times\Omega} \frac{|D^\beta v(x)-D^\beta v(y)|^2}{|x-y|^{d+2s}} \, \delta(x,y)^{2\kappa} \, dy \, dx
\]
and the associated space
\begin{equation} \label{eq:weighted_sobolev}
H^\ell_\kappa (\Omega) :=  \left\{ v \in H^\ell(\Omega) \colon \| v \|_{H^\ell_\kappa (\Omega)} < \infty \right\} .
\end{equation}

The regularity estimate in the weighted Sobolev scale \eqref{eq:weighted_sobolev} reads as follows.

\begin{theorem}[weighted Sobolev estimate] \label{T:weighted_regularity}
Let $\Omega$ be a bounded, Lipschitz domain satisfying the exterior ball condition, $s\in(0,1)$, $f \in C^{1-s}(\overline\Omega)$ and $u$ be the solution of \eqref{eq:weak_linear}. Then, for every $\eps>0$ we have $u \in \widetilde H^{s+1-2\eps}_{1/2-\eps}(\Omega)$ and
\[
\|u\|_{\widetilde H^{s+1-2\eps}_{1/2-\eps}(\Omega)} \le \frac{C(\Omega,d,s)}{\eps} \| f \|_{C^{1-s}(\overline\Omega)}.
\]
\end{theorem}

For simplicity, the theorem above was stated using the weight $\kappa = 1/2 -\eps$; a more general form of the result can be found also in \cite{BBNOS18}. We point out that, for its application in finite element analysis, such a choice is optimal. In principle, increasing the exponent $\kappa$ of the weight allows for a higher differentiability in the solution, with no restriction on $\kappa$ above (as long as $f$ is sufficiently smooth). However, when exploiting this weighted regularity by introducing approximations on a family of shape-regular and graded meshes, the order of convergence (with respect to the number of degrees of freedom) is only incremented as long as $\kappa < 1/2$.

It is worth pointing out that \Cref{T:weighted_regularity} is valid for Lipschitz domains satisfying an exterior ball condition. Although such a condition on the domain is much less restrictive than the $C^\infty$ requirement in \Cref{T:reg_grubb}, for polytopal domains it implies convexity. For that reason, we present here a result of an ongoing work \cite{BoNo19}, that characterizes regularity of solutions in terms of Besov norms.

\begin{theorem}[regularity on Lipschitz domains] \label{T:Besov_regularity}
Let $\Omega$ be a bounded Lipschitz domain, $s \in (0,1)$ and $f \in
H^r(\Omega)$ for some $r \in (-s,0]$. Then, the
solution $u$ to \eqref{eq:Dirichlet} belongs to the Besov space $B^{s+\gainB}_{2,\infty}(\Omega)$, where $\gainB = \min \{ s + r -\eps , 1/2 \}$, with $\eps > 0$ arbitrarily small, with
\begin{equation} \label{eq:Besov_regularity}
\| u \|_{B^{s+\gainB}_{2,\infty}(\Omega)} \le C(\Omega,d,s,\gainB) \|f\|_{H^r(\Omega)}.
\end{equation}
Consequently, by an elementary embedding, we deduce
\begin{equation} \label{eq:regularity}
\| u \|_{H^{s+\gain}(\Omega)} \le\frac{C(\Omega,d,s,\gain)}{\eps} \|f\|_{H^r(\Omega)},
\end{equation}
where $\gain = \min \{s+r, 1/2\} - \eps$ is `almost' as in \Cref{T:reg_grubb}.
\end{theorem}

We briefly outline the main ideas in the proof of \Cref{T:Besov_regularity}, which follows a technique proposed by Savar\'e \cite{Savare98} for local problems. The point is to use the classical Nirenberg difference quotient method, and thus bound a certain Besov seminorm of the solution $u$. 

Let $t \in (0,1)$ and $D$ be a set generating $\mRd$ and star-shaped with respect to the origin (for example, a cone). Then
the functional 
\begin{equation} \label{eq:Besov-norm}
  [v]_{s+t,2,\Omega} := \sup_{h \in D \setminus \{ 0 \}}
  \frac{1}{|h|^t} |v-v(\cdot+h)|_{H^s(\Omega)}
\end{equation}
induces the standard seminorm in the Besov space $B^{s+t}_{2,\infty}(\Omega)$. 
Because $\Omega$ is a Lipschitz domain, it satisfies a uniform cone property; upon a partition of unity argument, this gives (finitely many) suitable sets $D$ where translations can be taken. 
For a localized translation operator $T_h$, it is possible to prove a bound of the
form
\begin{equation} \label{eq:translation_bound}
|T_h u - u|_s^2 \le C \, |h|^s \, |u|_s^2,
\end{equation}
which, in view of \eqref{eq:Besov-norm}, yields $u \in B^{3s/2}_{2,\infty}(\Omega)$. Once this estimate has been obtained, a bootstrap argument leads to \eqref{eq:Besov_regularity}. Moreover, a refined estimate in $B^{3s/2}_{2,\infty}(\Omega)$ reads
\[
|u|_{B^{3s/2}_{2,\infty}(\Omega)} \lesssim \|f\|_{B^{-s/2}_{2,1}(\Omega)},
\]  
and interpolation with
$|u|_{\widetilde{H}^s(\Omega)} \lesssim \|f\|_{H^{-s}(\Omega)}$ yields the following result
\cite{BoNo19}.

\begin{theorem}[lift theorem on Lipschitz domains]\label{T:lift}
Let $\Omega$ be a bounded Lipschitz domain, $s \in (0,1)$ and $f \in
H^r(\Omega)$ for some $r \in (-s,-s/2]$. Then, the
solution $u$ to \eqref{eq:Dirichlet} belongs to the Sobolev space $\widetilde{H}^{r+2s}(\Omega)$, with
\[
\|u\|_{\widetilde{H}^{r+2s}(\Omega)} \lesssim \|f\|_{H^{r}(\Omega)}.
\]
\end{theorem}

\subsection{Finite element discretization} \label{sec:FE_linear}
In this section, we discuss a direct finite element method to approximate \eqref{eq:weak_linear} using piecewise linear continuous functions. We 
consider a family $\{\Th \}_{h>0}$ of conforming and simplicial meshes of $\Omega$, that we assume to be shape-regular, namely,
\[
  \sigma := \sup_{h>0} \max_{T \in \Th} \frac{h_T}{\rho_T} <\infty,
\]
where $h_T = \mbox{diam}(T)$ and $\rho_T $ is the diameter of the largest ball contained in $T$. As usual, the subindex $h$ denotes the mesh size, $h = \max_{T \in \Th} h_T$ and we take elements to be closed sets.
We denote by $\mathcal{N}_h$ the set of interior vertices of $\Th$, by $N$ the cardinality of $\mathcal{N}_h$, and by $\{ \varphi_i \}_{i=1}^N$ the standard piecewise linear Lagrangian basis, with $\phii_i$ associated to the node $\x_i \in \mathcal{N}_h$. Thus, the set of discrete functions is
\begin{equation} \label{eq:FE_space}
\mathbb{V}_h :=  \l\{ v \in C_0(\Omega) \colon v = \sum_{i=1}^N v_i \varphi_i \r\},
\end{equation}
and is clearly conforming: $\mathbb{V}_h \subset \tHs$ for all $s \in (0,1)$.

With the notation described above, the discrete counterpart to \eqref{eq:weak_linear} reads: find $u_h \in \mathbb{V}_h$ such that
\begin{equation} \label{eq:weak_linear_discrete}
(u_h, v_h)_s = \langle f, v_h \rangle \quad \forall v_h \in \mathbb{V}_h.
\end{equation}
Because $u_h$ is the projection of $u$ onto $\mathbb{V}_h$ with respect to the $\tHs$-norm, we have the best approximation property
\begin{equation} \label{eq:best_approximation_linear}
\|u - u_h \|_{\tHs} = \min_{v_h \in \mathbb{V}_h} \|u - v_h \|_{\tHs}.
\end{equation}

Therefore, in order to obtain a priori rates of convergence in the energy norm, it suffices to bound the distance between the discrete spaces and the solution. Although the bilinear form $(\cdot, \cdot)_s$ involves integration on $\Omega\times\mRd$, one can apply an argument based on a fractional Hardy  inequality, to prove that the energy norm may be bounded in terms of fractional--order norms on $\Omega$ (see \cite{AcosBort2017fractional}). It follows that bounding errors within $\Omega$ leads to error estimates in the energy norm.

A technical aspect of fractional-order seminorms is that they are not additive with respect to domain decompositions. With the goal of deriving interpolation estimates, we define the star or first ring of an element $T \in \Th$ by
\[
  S^1_T := \bigcup \left\{ T' \in \Th \colon \oT \cap \oTp \neq \emptyset \right\}.
\]
We also introduce the star of $S^1_T$ (or second ring of $T$),
\[
  S^2_T := \bigcup \left\{ T' \in \Th \colon S^1_T \cap \oTp \neq \emptyset \right\},
\]
and the star of the node $\x_i \in \mathcal{N}_h$, $S_i := \mbox{supp}(\varphi_i)$.

Faermann \cite{Faermann} proved the localization estimate 
\[
|v|_{H^s(\Omega)}^2 \leq \sum_{T \in \Th} \left[ \int_T \int_{S^1_T} \frac{|v (x) - v (y)|^2}{|x-y |^{d+2s}} \; dy \; dx + \frac{C(d,\sigma)}{s h_T^{2s}} \| v \|^2_{L^2(T)} \right] \quad \forall v \in H^s(\Omega). 
\]
This inequality shows that to estimate fractional seminorms over $\Omega$, it suffices to compute integrals over the set of patches $\{T \times S^1_T \}_{T \in \Th}$ plus local zero-order contributions. Bearing this in mind, one can prove the following type of estimates for suitable quasi-interpolation operators (see, for example, \cite{BoNoSa18,CiarletJr}).

\begin{proposition}[interpolation estimates on quasi-uniform meshes] \label{prop:app_SZ}
Let $T \in \T$, $s \in (0,1)$, $\ell \in (s, 2]$, and $\SZ$ be a suitable quasi-interpolation operator. 
  If $v \in H^\ell (S^2_T)$, then
\begin{equation} \label{eq:interpolation}
 \int_T \int_{S^1_T} \frac{|(v-\SZ v) (x) - (v-\SZ v) (y)|^2}{|x-y|^{d+2s}} \, d y \, d x \le C \, h_T^{2(\ell-s)} |v|_{H^\ell(S^2_T)}^2.
\end{equation}
where $C = C(\Omega,d,s,\sigma, \ell)$.
Therefore, for all $v \in H^\ell (\Omega)$, it holds
\begin{equation} \label{eq:global_interpolation}
| v - \SZ v|_{\Hs} \le C(\Omega,d,s,\sigma, \ell) \, h^{\ell-s} |v|_{H^\ell(\Omega)}.
\end{equation}
\end{proposition}

The statement \eqref{eq:global_interpolation} in \Cref{prop:app_SZ} could have also been obtained by interpolation of standard integer-order interpolation estimates. However, the technique of summing localized fractional-order estimates also works for graded meshes (cf. \eqref{eq:weighted_interpolation} and \eqref{eq:global_weighted_interpolation} below).

Combining estimate \eqref{eq:global_interpolation} with the best approximation property \eqref{eq:best_approximation_linear} and the regularity estimates described in \Cref{sec:regularity_linear}, we can derive convergence rates. Concretely, the estimates from \Cref{T:reg_grubb} and \Cref{T:Besov_regularity} translate into a priori rates for quasi-uniform meshes. However, optimal application of \Cref{T:weighted_regularity} requires a certain type of mesh grading.
In two-dimensional problems ($d=2$), this can be attained by constructing \emph{graded} meshes in the spirit of \cite[Section 8.4]{Grisvard}. In addition to shape regularity, we assume that the family $\{\Th \}$ satisfies the following property: there is a number $\mu\ge1$ such that given a parameter $h$ representing the meshsize at distance $1$ to the boundary $\partial\Omega$ and $T\in\Th$, we have
 \begin{equation} \label{eq:H}
 h_T \leq C(\sigma)
 \begin{dcases}
   h^\mu, & T \cap \partial \Omega \neq \emptyset, \\
   h \dist(T,\pp \Omega)^{(\mu-1)/\mu}, & T \cap \partial \Omega = \emptyset.
 \end{dcases}
\end{equation} 
The number of degrees of freedom is related to $h$ by means of the parameter $\mu$ because (recall that $d=2$)
\[
 N = \dim \mathbb{V}_h \approx
  \begin{dcases}
    h^{-2}, & \mu < 2, \\
    h^{-2} | \log h |, & \mu = 2, \\
    h^{-\mu}, & \mu > 2.
  \end{dcases}
\]
Also, $\mu$ needs to be related to the exponent $\kappa$ used in the Sobolev regularity estimate (cf. \Cref{T:weighted_regularity}). It can be shown that the choice $\mu = 2$, that corresponds to $\kappa = 1/2$, yields optimal convergence rates in terms of the dimension of $\mathbb{V}_h$. We also remark that, as discussed in \cite{BoCi19}, for three-dimensional problems ($d=3$), the grading strategy \eqref{eq:H} becomes less flexible, and yields lower convergence rates. For optimal mesh grading beyond $\mu=2$ for both $d=2,3$ we need to break the shape regularity assumption and resort to anisotropic finite elements. They in turn are less flexible in dealing with the isotropic fractional norm of $H^s(\Omega)$ and its localization \cite{Faermann}. This important topic remains open.

Quasi-interpolation estimates in weighted Sobolev spaces \eqref{eq:weighted_sobolev} can be derived in the same way as in \Cref{prop:app_SZ}. More precisely, the weighted counterparts to \eqref{eq:interpolation} and \eqref{eq:global_interpolation} read 
\begin{equation} \label{eq:weighted_interpolation}
 \int_T \int_{S^1_T} \frac{|(v-\SZ v) (x) - (v-\SZ v) (y)|^2}{|x-y|^{d+2s}} \, d y \, d x \le C
   h_T^{2(\ell-s-\kappa)} |v|_{H^\ell_\kappa(S^2_T)}^2,
\end{equation}
for all $v \in H^\ell_\kappa (S^2_T)$ and
\begin{equation} \label{eq:global_weighted_interpolation}
| v - \SZ v|_{\Hs} \le C  h^{\ell-s-\alpha} |v|_{H^\ell_\kappa(\Omega)} \quad \forall v \in H^\ell_\kappa (\Omega),
\end{equation}
respectively. The constants above are $C = C(\Omega,d,s,\sigma,\ell,\kappa)$. 

We collect all the convergence estimates in the energy norm --involving quasi-uniform and graded meshes-- in a single statement \cite{AcosBort2017fractional}.

\begin{theorem}[energy error estimates for linear problem] \label{T:conv_linear}
Let $u$ denote the solution to \eqref{eq:weak_linear} and denote by $u_h \in \mathbb{V}_h$ the solution of the discrete problem \eqref{eq:weak_linear_discrete}, computed over a mesh $\T$ consisting of elements with maximum diameter $h$. If $f \in L^2(\Omega)$, we have
\[
\|u - u_h \|_{\tHs}  \le C(\Omega,d,s,\sigma) \,  h^\alpha |\log h| \, \|f\|_{L^2(\Omega)},
\]
where $\alpha = \min \{s, 1/2 \}$.
Additionally, if $\Omega \subset \mR^2$, $f \in C^{1-s}(\overline{\Omega})$ and the family $\{\Th\}$ satisfies \eqref{eq:H} with $\mu = 2$, we have
\[
\|u - u_h \|_{\tHs}  \le C(\Omega,s,\sigma) \, h |\log h| \|f\|_{C^{1-s}(\overline{\Omega})}.
\]
\end{theorem}

To illustrate that \Cref{T:conv_linear} is sharp, we solve the problem from \Cref{ex:nonsmooth} on the discrete spaces \eqref{eq:FE_space} using a family of quasi-uniform meshes and a family of meshes graded according to \eqref{eq:H}. In \Cref{tab:ejemplo}, we report computational convergence rates in the energy norm for several values of $s$. 
We observe good agreement with the rates predicted by \Cref{T:conv_linear}.
\begin{table}[htbp]\small\centering
\begin{tabular}{|c| c| c| c| c| c| c| c| c| c|} \hline
Value of $s$ & 0.1 & 0.2 & 0.3 & 0.4 & 0.5 & 0.6 & 0.7 & 0.8 & 0.9 \\ \hline
Uniform meshes & 0.497 & 0.496 &  0.498  & 0.500 & 0.501  & 0.505 & 0.504  & 0.503 &  0.532 \\ \hline
Graded meshes & 1.066 & 1.040  & 1.019 & 1.002  & 1.066 & 1.051  & 0.990 & 0.985 & 0.977  \\  \hline   
\end{tabular}
\bigskip
\caption{Computational rates of convergence (with respect to $h$) for the problem from \Cref{ex:nonsmooth} in $d=2$ dimensions. Rates using quasi-uniform meshes are listed in the second row, while rates using graded meshes, with $\mu = 2$ in \eqref{eq:H}, are reported in the third row.}
\label{tab:ejemplo}
\end{table}

\subsection{Computational challenges} \label{sec:computational_linear}
Having at hand theoretical estimates for finite element discretizations of \eqref{eq:weak_linear}, we still need to address how to compute discrete solutions and, in particular, how to accelerate the assembly and solution of the discrete system that arises.

\medskip\noindent
{\bf Matrix assembly.}
We first comment on key aspects of the finite element implementation for problems in dimension $d = 2$. If $\vU = (u_i)_{i=1}^N$ and $u_h = \sum_{i=1}^N u_i \varphi_i$, it follows from \eqref{eq:FE_space} and \eqref{eq:weak_linear_discrete} that the linear finite element system reads $\vA \vU = \vF$ with stiffness matrix $\vA$ and right-hand side vector $\vF$ given by
\[
\vA_{ij} = (\phii_i, \phii_j)_s, \quad \vF_i = \langle f, \phii_i \rangle.
\]

Computation of the stiffness matrix is not an easy task. There are two numerical difficulties in taking a direct approach. In first place, the bilinear form $(\cdot, \cdot)_s$ requires integration on unbounded domains; we point out that --at least for homogeneous problems as the ones considered here-- integration over $\Omega\times\Omega^c$ can be reduced to a suitable integration over $\Omega\times\pp\Omega$ by using the Divergence Theorem \cite{AiGl17}. Secondly, suitable quadrature rules to compute the stiffness matrix entries are required. To handle the singular (non-integrable) kernel $|x|^{-d-2s}$, one could use techniques from the boundary element method \cite{ChernovPetersdorffSchwab:11,SauterSchwab}; we refer to \cite{AcosBersBort2017short} for details.

\medskip\noindent
{\bf Compression.}
Note that, independently of $s$, finite element spaces \eqref{eq:FE_space} give rise to {\em full} stiffness matrices. Indeed, for any pair of nodes $\x_i, \x_j$ such that $S_i \cap S_j = \emptyset$, 
\[
\vA_{ij} =  - C_{d,s} \iint_{S_i \times S_j} \frac{\phii_i(x) \; \phii_j(y)}{|x-y|^{d+2s}} \; dy \; dx < 0.
\]
Thus, computation of the stiffness matrix $\vA$ involves a large number of far-field interactions, that is, elements $\vA_{ij}$ for $\x_i$ and $\x_j$ sufficiently far. However, these elements should be significantly smaller than the ones that involve nodes close to one another. In \cite{AiGl17, zhao2017adaptive} the cluster method from the boundary element literature was applied and, instead of computing and storing all individual elements from $\vA$, far field contributions are replaced by suitable low-rank blocks. The resulting data-sparse representation has $\mathcal{O}(N \log^\alpha N)$ complexity for some $\alpha \ge 0$. Reference \cite{karkulik2018mathcal} shows that the inverse of $\vA$ can be represented using the same block structure as employed to compress the stiffness matrix.

\medskip\noindent
{\bf Preconditioning.}
There are also issues to be addressed regarding the solution of the dense matrix equation $\vA \vU = \vF$. The use of matrix factorization to solve such a system has complexity $\mathcal{O}(N^3)$. As an alternative, one can use a conjugate gradient method and thereby the number of iterations needed for a fixed tolerance scales like $\sqrt{\kappa(\vA)}$, where $\kappa(\vA)$ is the condition number of $\vA$ and satisfies \cite{AiMcTr99}
\[
\kappa(\vA) = \mathcal{O}\left( N^{2s/d} \left( \frac{h_{max}}{h_{min}} \right)^{d-2s} \right). 
\]
Therefore, for two-dimensional problems, we deduce $\kappa(\vA) = \mathcal{O}(h^{-2s})$ for quasi-uniform meshes, while $\kappa(\vA) = \mathcal{O}(h^{-2} |\log h|^s)$ for meshes graded according to \eqref{eq:H} with $\mu = 2$. In the latter case, diagonal preconditioning allows us to recover the same condition number as for uniform meshes \cite{AiMcTr99}.

Recently, there have been some advances in the development of preconditioners for fractional diffusion. For instance, multigrid preconditioners were mentioned in \cite{AiGl17}, while operator preconditioners were studied in \cite{gimperlein2019optimal}. We now briefly comment on some features of an additive Schwarz preconditioner of BPX-type \cite{fractionalbpx} (see also \cite{faustmann2019optimal}). Assume we have a hierarchy of discrete spaces $\mathbb{V}_0 \subset \ldots \mathbb{V}_J = \mathbb{V}$, with mesh size $h_j = \gamma^{2j}$, and let $\iota_j : \mathbb{V}_j \to \mathbb{V}$ be the inclusion operator. The basic ingredients needed to apply the general theory for additive Schwarz preconditioners are:
\begin{itemize}[leftmargin=*]
\item {\it Stable decomposition}: for every $v\in \mathbb{V}$, there exists a decomposition $v = \sum_{j=0}^J v_j$ with $v_j \in \mathbb{V}_j$ such that
\begin{equation} \label{eq:stable-decomposition}
\sum_{j=0}^J h_j^{-2s} \|v_j\|_{L^2(\Omega)}^2 \leq c_0 \|v\|_{\tHs}^2.
\end{equation}
A fundamental ingredient to prove this estimate for polyhedral domains is the optimal regularity pickup estimate for Lipschitz domains of \Cref{T:lift}, that allows us to perform an Aubin-Nitsche duality argument.

\item {\it Boundedness}: for every $v = \sum_{j=0}^J v_j$ with $v_j \in \mathbb{V}_j$,
\begin{equation} \label{eq:boundedness} 
\| \sum_{j=0}^J v_j \|_{\tHs}^2 \leq c_1 \sum_{j=0}^J h_j^{-2s} \|v_j\|_{L^2(\Omega)}^2.
\end{equation}
As usual, boundedness of multilevel decompositions can be proved by estimating how much scales interact (i.e., using a strengthened Cauchy-Schwarz inequality). Nonlocality adds some difficulties to the derivation of such an estimate, because one cannot integrate by parts elementwise. The argument in \cite{fractionalbpx} is based on the Fourier representation of the fractional Laplacian.
\end{itemize}
The conditions \eqref{eq:stable-decomposition} and \eqref{eq:boundedness} imply that the preconditioner $\vB := \sum_{j=0}^J h_j^{2s-d} \iota_j \iota'_j$ satisfies $\kappa(\vB \vA) \le \frac{c_0}{c_1}$ for graded bisection grids \cite{fractionalbpx}. We
  illustrate this statement in Table \ref{tab:FL-BPX-bisect} for $\Omega=(-1,1)^2$,
  $f=1$ and $s=0.9, 0.5, 0.1$. We observe a mild increase of iteration counts but
  rather robust performance with respect to $s$. 
\begin{table}[!ht]
\centering
\begin{tabular}{|m{0.4cm}|m{1.cm}|m{0.75cm}|m{0.75cm}|m{0.75cm}
  |m{0.75cm}|m{0.75cm}|m{0.75cm}|m{0.75cm}|m{0.75cm}|m{0.75cm}|}
  \hline
  \multirow{2}{*}{$\bar{J}$} & \multirow{2}{*}{$N$} 
  & \multicolumn{3}{c|}{$s=0.9$}
  & \multicolumn{3}{c|}{$s=0.5$} 
  & \multicolumn{3}{c|}{$s=0.1$}\\ \cline{3-11}
  & & GS & CG  & PCG  & GS & CG & PCG & GS & CG & PCG \\ \hline
  8 &  209 & 101 & 18 & 17 & 22 & 14 & 16 & 9 & 21 & 18 \\ \hline
  9 &  413 & 172 & 24 & 20 & 30 & 18 & 18 & 9 & 25 & 21 \\ \hline
  10 & 821 & 314 & 32 & 22 & 42 & 20 & 19 & 9 & 29 & 23 \\ \hline
  11 & 1357 & 494 & 41 & 23 & 55 & 24 & 20 & 9 & 31 & 24 \\ \hline
  12 & 2753 & 792 & 55 & 25 & 72 & 30 & 21 & 9 & 36 & 26 \\ \hline
  13 & 4977 & 1391 & 73 & 26 & 94 & 34 & 22 & 10 & 37 & 27 \\ \hline
  14 & 9417 & 2357 & 95 & 27 & 133 & 39 & 23 & 10 & 40 & 28 \\ \hline
\end{tabular}
\bigskip
\caption{Number of iterations for Gauss-Seidel (GS), conjugate gradient (CG)
  and preconditioned CG with BPX preconditioner (PCG). Stopping criteria is
    $\frac{\|\vA\vU - \vF\|_2}{\|\vF\|_2} < 10^{-6}$.}
\label{tab:FL-BPX-bisect}
\end{table}
%

\subsection{Applications and related problems} \label{sec:applications_linear}
Numerical methods for fractional diffusion models have been extensively studied recently. Let us mention some applications on linear elliptic problems of the approach treated in this section:

\begin{itemize}[leftmargin=*]
\item {\em A posteriori error analyisis and adaptivity:} The reduced regularity of solutions of \eqref{eq:weak_linear} and the high computational cost of assembling the stiffness matrix $\vA$ motivate the pursue of suitable adaptive finite element methods. A posteriori error estimates of residual type have been proposed and analyzed in \cite{AiGl17,faustmann2019quasi,gimperlein2019space,NvPZ}, and gradient-recovery based estimates in \cite{zhao2017adaptive}.
Adaptivity is however a topic of current research \cite{faustmann2019optimal,faustmann2019quasi,gimperlein2019space}.

\item {\em Eigenvalue problems:} The fractional eigenvalue problem arises, for example, in quantum mechanics problems in which the Brownian-like quantum paths are replaced by L\'evy-like ones in the Feynman path integral \cite{Laskin}. Reference \cite{BdPM} studies conforming finite element approximations and applies the Babu\v{s}ka-Osborn theory \cite{BO91}, thereby obtaining convergence rates for eigenfunctions (in the energy and in the $L^2$-norms) and eigenvalues. Other methods, implemented on one-dimensional problems, include finite differences \cite{duo2015computing} and matrix methods \cite{ghelardoni2017matrix, zoia2007fractional}.

\item {\em Control problems:} Finite element methods for linear-quadratic optimal control problems involving the fractional Laplacian \eqref{eq:def_Laps} have been studied recently. In these problems, the control may be located inside \cite{d2018priori} or outside the domain \cite{AntilKhatriWarma18}.

\item {\em Non-homogeneous Dirichlet conditions:} 
A mixed method for the non-homogeneous Dirichlet problem for the integral fractional Laplacian was proposed in \cite{AcBoHe18}. Such a method is based on weak enforcement of the Dirichlet condition and using a suitable non-local derivative \cite{DROV17} as a Lagrange multiplier. To circumvent approximating the nonlocal derivative, \cite{AntilKhatriWarma18} proposed approximations of the non-homogeneous Dirichlet problem by a suitable Robin exterior value problem.
\end{itemize}

\subsection{Nonconforming FEM: Dunford-Taylor approach} \label{sec:dunford-taylor}
We finally report on a finite element approach for \eqref{eq:Dirichlet} proposed in \cite{BoLePa17} and based on the Fourier representation of the $\tHs$-inner product:
\[
(v,w)_s = \int_{\mRd} |\xi|^s \Fou(v) |\xi|^s \overline{\Fou(w)} d\xi
= \int_{\mRd} \Fou((-\Delta)^s v)(\xi) \overline{\Fou(w(\xi))} d\x.
\]
This expression can be equivalently written as
\begin{equation} \label{eq:DT}
(v,w)_s = \frac{2\sin(s\pi)}{\pi}
\int_0^\infty t^{1-2s} \int_{\R^d} \big( -\Delta
(I-t^2\Delta)^{-1} v  \big) w \, dx dt.
\end{equation}
To see this, use Parseval's formula to obtain
\[
\int_{\R^d} \big( -\Delta
(I-t^2\Delta)^{-1} v  \big) w \, dx  = \int_{\R^d}
\frac{|\xi|^2}{1+t^2|\xi|^2} \Fou(v)(\xi) \overline{\Fou(w)(\xi)} d\xi,
\]
followed by the change of variables $z=t|\xi|$, which converts the repeated integrals
in the expression for $(v,w)_s$ into separate integrals, one of them being
\[
\int_0^\infty \frac{z^{1-2s}}{1+z^2} dz = \frac{\pi}{2\sin(s \pi)}.
\]
Although identity \eqref{eq:DT} is not an integral representation of the operator 
$\Laps$, but rather of the bilinear form $(\cdot, \cdot)_s$, we regard it as a Dunford-Taylor representation.
  
To set up this formal calculation in the correct functional framework, given
$u\in \tHs \subset L^2(\R^d)$ and $t>0$, let  $v(u,t) \in H^{2+s}(\R^d)$
be the solution to 
$v - t^2 \Delta v = -u$ in $\mRd$, or equivalently $v = -(I-t^2\Delta)^{-1}u$.
Therefore, $\Delta v = t^{-2} (v+u)$ and 
\[
(u,w)_s = \frac{2\sin(s \pi)}{\pi} \int_0^\infty t^{-1-2s} \langle u + v(u,t) , w \rangle \; dt \quad \forall u,w \in \tHs.
\]

This representation is the starting point of a three-step numerical method
\cite{BoLePa17}.
\begin{itemize}[leftmargin=*]
\item {\em Sinc quadrature:} the change of variables $t = e^{-y/2}$ yields
\[ (u,w)_s = \frac{\sin(s\pi)}{\pi} \int_{-\infty}^\infty e^{sy}  \langle u + v(u,t(y)) , w \rangle \; dy
\]
Thus, given an integer $N>0$ and a set of points $\{y_j\}_{j=-N}^N$ with uniform
  spacing $\approx N^{-1}$,
the sinc quadrature $Q_s(u,w)$ approximation of $(\cdot,\cdot)_s$ is given by
\[
Q_s(u,w) = \frac{\sin(s\pi)}{N \pi} \sum_{j = -N}^{N} e^{sy_j}
\langle u + v(u,t(y_j)), w \rangle.
\]

\item {\em Domain truncation:} We stress that, in spite of $u$ being supported in $\Omega$, $v(u,t)$ is supported in all of $\mRd$ for all $t$; hence, some truncation is required. The method from \cite{BoLePa17} considers, for given $M>0$,
a family of balls $B^M(t)$ that contain $\Omega$ and whose radius depends on $M$ and $t$ and can be computed a priori.

\item {\em Finite element approximation:} Finally, a standard finite element discretization on $B^M(t)$ is performed. This requires meshes that fit $\Omega$ and $B^M(t) \setminus \Omega$ exactly -- a non-trivial task; let us denote the discrete spaces on $\Omega$ and $B^M(t)$ by $\mathbb{V}_h$ and $\mathbb{V}_h^M$, respectively. Given $\psi \in L^2(\mRd)$, $t>0$ and $M>0$, we define $v_h^M = v_h^M(\psi, t)\in\mathbb{V}_h^M$ to be the unique solution of
\[
\int_{B^M(t)} v_h^M w_h + t^2 \nabla v_h^M \cdot \nabla w_h  \; dx = - \int_{B^M(t)} \psi w_h  \; dx \quad \forall \, w_h \in \mathbb{V}_h^M.
\]
\end{itemize}

The fully discrete bilinear form reads:
\[
a_\T^{N,M}(u_h, w_h) := \frac{\sin(s\pi)}{N \pi} \sum_{j =-
  N}^{N} e^{sy_j} \langle 
  u_h + v_h^M(u_h,t(y_j)), w_h
  \rangle \quad\forall \, u_h, w_h \in \mathbb{V}_h.
\]
Using a Strang's type argument to quantify the consistency errors generated
by the three steps above, one obtains the a priori estimate
\cite[Theorem 7.7]{BoLePa17}
\[
\| u - u_h \|_{\tHs} \le C \l( e^{-c\sqrt{N}} + e^{-cM} + h^{\beta-s} |\log h| \r) \| u \|_{\tH^\beta(\Omega)}, 
\]
where $\beta \in (s,3/2)$. Choosing $\beta = s + 1/2 -\eps$, which is consistent with \Cref{T:reg_grubb} and \Cref{T:Besov_regularity}, $M = \mathcal{O}(|\log h|)$ and $N = \mathcal{O}(|\log h|^2)$ gives the convergence rate
\[
\| u - u_h \|_{\tHs} \le C h^{\min\{s,\frac12\}} |\log h| \, \|f\|_{L_2(\Omega)}.
\]
This is similar to the rate obtained in \Cref{T:conv_linear} for quasi-uniform meshes. To the best of the authors' knowledge, implementation of this 
approach over graded meshes, while feasible in theory, has not yet been pursued in practice.

%% file: obstacle.tex
\section{Fractional Obstacle Problem} \label{sec:obstacle}
In this section we review finite element methods for the solution of the obstacle problem for the integral fractional Laplacian which, from now on, we shall simply refer to as the fractional obstacle problem.

The fractional obstacle problem appears, for example, in optimal stopping times for jump processes. In particular, it is used in the modeling of the rational price of perpetual American options \cite{ContTankov}. More precisely, if $u$ represents the rational price of a perpetual American option where the assets prices are modeled  by a L\'evy process $X_t$, if $\chi$ denotes the payoff function, then $u$ solves a fractional obstacle problem with obstacle $\chi$.

An a posteriori error analysis of approximations of variational inequalities involving integral operators on arbitrary bounded domains was performed in \cite{NvPZ}. We also comment on two recent works related to the approach we review here. Reference \cite{BG18} deals with finite element discretizations to obstacle problems involving finite and infinite-horizon nonlocal operators. The experiments shown therein were performed on one-dimensional problems with uniform meshes, and indicate convergence with order $h^{1/2}$ in the energy norm. A theoretical proof of that convergence order was obtained in \cite{BoLeSa18}, where approximations using the 
approach discussed in \Cref{sec:dunford-taylor} were considered. We also refer to \cite{gimperlein2019space} for computational comparisons between adaptive strategies and uniform and graded discretizations in two-dimensional problems.

\subsection{Variational formulation} \label{sec:formulation_obstacle}
As before, we assume that $\Omega \subset \mRd$ is an open and bounded domain and, for the sake of applying weighted regularity estimates, we assume that $\Omega$ has a Lipschitz boundary and satisfies the exterior ball condition. 
Given $s \in (0,1)$ and functions $f: \Omega \to \mR$ and $\chi : \overline\Omega \to \mR$, with $\chi < 0$ on $\pp\Omega$, the obstacle problem is a constrained minimization problem on $\tHs$ associated with a quadratic functional. Defining the admissible convex set
\[
\K := \left\{ v \in \tHs: v \geq \chi \mbox{ a.e in } \Omega \right\},
\]
the solution to the fractional obstacle problem is $u=\argmin_{v\in\K} \mathcal{J}(v)$
where
\[
\mathcal{J}(v) := \frac12 \| v \|_{\tHs}^2 - \langle f, v \rangle.
\]
Existence and uniqueness of solutions is standard. Taking first variation of $\mathcal{J}$, we deduce that such a minimizer $u \in \K$ solves the variational inequality
\begin{equation}
\label{eq:obstacle}
  (u,u-v)_s \leq \langle f, u-v \rangle \quad \forall v \in \K.
\end{equation}

It can be shown \cite{musina2017variational} that, if $f \in L^p(\Omega)$ for $p > d/2s$, then the solution to the obstacle problem is indeed a continuous function, and that it satisfies the complementarity condition
\begin{equation}
\label{eq:complementarity}
  \min\left\{ \lambda, u-\chi \right\} = 0 \quad\mbox{ a.e in } \Omega,  \mbox{ where } \quad \lambda := \Laps u - f.
\end{equation}
For our discussion, we assume that $f$ is such that the solution is defined pointwise, and consequently,  we define the coincidence (or contact) and non-coincidence sets,
\[
\Lambda := \{ x \in \Omega : u(x) = \chi(x) \}, \quad N := \Omega \setminus \Lambda.
\]
The complementarity condition \eqref{eq:complementarity} can be succinctly expressed as $\lambda \ge 0$ in $\Lambda$ and $\lambda = 0$ in $N$. The set $\pp \Lambda$, where the solution detaches from the obstacle, is the {\it free boundary.}

\subsection{Regularity} \label{sec:regularity_obstacle}
The following regularity results for solutions to the fractional obstacle problem are instrumental for error analysis. We recall our assumption that the obstacle $\chi$ is a continuous function and strictly negative on $\pp\Omega$:
\begin{equation} \label{eq:cond_chi}
  \varrho := \dist\left( \{ \chi>0\}, \pp \Omega \right) > 0.
\end{equation}
Furthermore, we shall assume that $f \ge 0$. Heuristically, these assumptions should guarantee that the behavior of solutions near $\pp\Omega$ is dictated by a linear problem and that the nonlinearity is confined to the interior of the domain. Finally, to derive regularity estimates, we assume that the data satisfy
\begin{equation}
\label{eq:defofcalF}
\chi \in C^{2,1}(\Omega), \quad f \in  \F(\overline\Omega) = \begin{dcases}
         C^{2,1-2s+\epsilon}(\overline\Omega), & s\in \left(0,\frac12\right], \\
         C^{1,2-2s+\epsilon}(\overline\Omega), & s \in \left(\frac12,1\right),
       \end{dcases}
\end{equation}
where $\epsilon>0$ is sufficiently small, so that $1-2s+\epsilon$ is not an integer.
Under these conditions, Caffarelli, Salsa and Silvestre \cite{CaSaSi08} proved that the solution to the problem posed in the whole space (with suitable decay conditions at infinity) is of class $C^{1,s}(\R^d)$. It is worth examining the limiting cases $s=1$ and $s=0$. The former corresponds to the classical obstacle problem whose solutions are of class $C^{1,1}(\R^d)$. The latter reduces to $\min\{u-\chi,u-f\}=0$ whose solutions are just of class $C^{0,1}(\R^d)$. The regularity of \cite{CaSaSi08} is thus a natural intermediate result.

We emphasize that deriving interior regularity estimates for \eqref{eq:obstacle} from this result, which is valid for problems posed in $\mRd$, is not as straightforward as for classical problems. Indeed, the nonlocal structure of $\Laps$ implies that, if $0\le\eta\le1$ is a smooth cut-off function such that $\eta=1$ in $\{ \chi > 0 \}$, then
\[
  \Laps (\eta u) \ne \eta \Laps u \quad \mbox{in } \{ \eta = 1 \}.
\]
To overcome this difficulty, reference \cite{BoNoSa18} proceeds as follows. Given a set $D$ such that $\{\chi>0\}\subset D \subset \Omega$, one can define a cutoff $\eta$ such that $D \subset \{ \eta = 1 \}$ and  split the space roughly into a region where $\eta = 1$, a region where $\eta = 0$ and a transition region. In the first two regions, $\Laps(\eta u)$ essentially coincides with a convolution operator with kernel $|z|^{-d-2s}$ but regularized at the origin, while the latter region is contained in the non-contact set $N$ and allows one to invoke interior regularity estimates for linear problems involving $\Laps$. An important outcome is that solutions to fractional obstacle problems are more regular near the free boundary ($C^{1,s}$) than near the domain boundary ($C^{0,s}$). This is critical for approximation.

Alternatively, one may invoke the Caffarelli-Silvestre extension \cite{CS:07} to obtain local regularity estimates \cite{CaSaSi08}. Since the extension problem involves a degenerate elliptic equation with a Muckenhoupt weight of class $A_2$ that depends only on the extended variable, one needs to combine fine estimates for degenerate equations with the translation invariance in the $x$-variable of the Caffarelli-Silvestre weight.

Once the interior regularity of solutions is established, one can invoke the H\"older boundary estimates for linear problems \cite{RosOtonSerra} and perform an argument similar to the one in \cite{AcosBort2017fractional} to deduce weighted Sobolev regularity estimates \cite{BoNoSa18}.

\begin{theorem}[weighted Sobolev regularity for the obstacle problem]
\label{T:regularity_obstacle}
Let $\Omega$ be a bounded Lipschitz domain satisfying the exterior ball condition, $s \in (0,1)$, and $\chi \in C^{2,1}(\Omega)$ satisfying \eqref{eq:cond_chi}.
Moreover, let $0 \leq f \in \F(\overline\Omega)$ and $u \in \tHs$ be the solution to \eqref{eq:obstacle}. For every $\eps >0$ we have that $u \in \tH^{s+1-2\eps}_{1/2-\eps}(\Omega)$ with the estimate
\begin{equation} \label{eq:regularity_obstacle}
  \|u\|_{\tH^{s+1-2\eps}_{1/2-\eps}(\Omega)} \leq \frac{C(\chi, s, d, \Omega, \varrho, \| f \|_{\F(\overline\Omega)})}\eps.
\end{equation}
\end{theorem}

We have stated the estimate in \Cref{T:regularity_obstacle} in weighted spaces because we are interested in the application of that result for finite element schemes over graded meshes. With the same arguments as in \cite{BoNoSa18}, it can be shown that the solution to the fractional obstacle problem \eqref{eq:obstacle} satisfies $u \in \tH^{s+1/2-\eps}(\Omega)$ and
\[
  \|u\|_{\tH^{s+1/2-\eps}(\Omega)} \leq \frac{C(\chi, s, d, \Omega, \varrho, \| f \|_{\F(\overline\Omega)})}\eps.
\]

A similar result, for the obstacle problem for a class of integro-differential operators, was obtained in \cite{BoLeSa18}. In the case of purely fractional diffusion (i.e., problems without a second-order differential operator), the estimate builds on \cite{Grubb} (cf. \Cref{T:reg_grubb}). Therefore, we point out that using \Cref{T:Besov_regularity}, the requirement that $\Omega$ be a $C^\infty$ domain in \cite[Cases A and B]{BoLeSa18} can be relaxed to $\Omega$ being Lipschitz.

\subsection{Finite element discretization} \label{sec:FE_obstacle}
We consider the same finite element setting as in \Cref{sec:FE_linear}: let $\mathbb{V}_h$ be linear Lagrangian finite element spaces as in \eqref{eq:FE_space} over a family of conforming and simplicial meshes $\T$. An instrumental tool in the analysis we review here is the interpolation operator $\SZ : L^1(\Omega) \to \mathbb{V}_h$ introduced in \cite{ChenNochetto}  that, besides satisfying \eqref{eq:interpolation}, is {\it positivity preserving}: it satisfies $\SZ v \ge 0$ for all $v \ge 0$. Such a property yields that, for every $v \in \K$,
\begin{equation} \label{eq:admissible_projection}
\SZ v \ge \SZ \chi \ \mbox{ in } \Omega.
\end{equation}

We therefore define the discrete admissible convex set 
\[
\K_h := \left\{ v_h \in \mathbb{V}_h: v_h \geq \SZ \chi \mbox{ in } \Omega \right\},
\]
and consider the discrete fractional obstacle problem: find $u_h \in \K_h$ such that 
\begin{equation}
\label{eq:obstacle_discrete}
  (u_h,u_h-v_h)_s \leq \langle f, u_h-v_h \rangle \quad \forall v_h \in \K_h.
\end{equation}
 
We illustrate the delicate interplay between regularity and approximability next. We exploit that $u$ is both globally of class $C^{0,s}(\overline{\Omega})$, via graded meshes as in the linear problem, and locally of class $C^{1,s}(\Omega)$. 
First, we split the error as
$\| u-u_h \|_{\tHs}^2 = (u - u_h, u - I_h u)_s + (u - u_h, I_h u - u_h)_s$, use Cauchy-Schwarz inequality and the interpolation estimate \eqref{eq:global_weighted_interpolation} to deduce
\[
\frac12 \|u-u_h \|_{\tHs}^2 \le 
C h^{2(1-2\eps)} \|u \|^2_{\widetilde H^{1+s-2\eps}_{1/2-\eps}(\Omega)} + (u - u_h, I_h u - u_h)_s.
\]
This is a consequence of Theorem \ref{T:weighted_regularity} and the use of graded meshes with parameter $\mu=2$ as in the linear theory of Section \ref{sec:FE_linear}. For the remaining term we integrate by parts and utilize the discrete variational inequality \eqref{eq:obstacle_discrete} to arrive at
\[
\begin{aligned}
( u - u_h &, I_h u - u_h )_s \le \int_\Omega (I_h u - u_h) \big((-\Delta)^su -f \big) \\
 & = \int_\Omega \Big[ (u -\chi) + \underbrace{(I_h \chi - u_h)}_{\le 0} + \big(I_h (u - \chi) - (u -\chi) \big) \Big] \Big(\underbrace{(-\Delta)^su -f}_{\ge 0} \Big).
\end{aligned}
\]
Invoking the complementarity condition \eqref{eq:complementarity}, we obtain
\[
(u - u_h, I_h u - u_h)_s \le \sum_{T \in \T}  \int_T \big( I_h (u - \chi) - (u - \chi) \big)  \big((-\Delta)^su -f \big).
\]
We next observe that the integrand does not vanish only for elements $T$ in the vicinity of the free boundary, namely $T$'s for which $u\ne\chi$ and $(-\Delta)^su \ne f$. Exploiting that $u\in C^{1,s}(\Omega)$, we infer that $(-\Delta)^su -f \in C^{0,1-s}(\Omega)$, whence
\[
\big| \big( \Laps u  -f \big) \big( I_h (u - \chi) -
(u - \chi) \big) \big| \le C h^{2} .
\]
This yields the following optimal energy error estimate. We refer to \cite{BoNoSa18} for details.

\begin{theorem}[error estimate for obstacle problem]
\label{thm:conv_rates}
Let $u$ be the solution to \eqref{eq:obstacle} and $u_h$ be the solution to \eqref{eq:obstacle_discrete}, respectively. Assume that $\chi \in C^{2,1}(\Omega)$ satisfies \eqref{eq:cond_chi} and that $f \in \F(\overline\Omega)$. If $d=2$, $\Omega$ is a convex polygon, and the meshes satisfy the grading hypothesis \eqref{eq:H} with $\mu = 2$, then we have that
\[ \begin{aligned}
 & \|u-u_h\|_{\tHs} \leq C h |\log h| & (s \ne 1/2), \\
 & \|u-u_h\|_{\widetilde H^{1/2}(\Omega)} \leq C h |\log h|^2 & (s = 1/2),
\end{aligned} \]
where $C>0$ depends on $\chi$, $s$, $d$, $\Omega$, $\varrho$ and $\| f \|_{\F(\overline\Omega)}$.
\end{theorem}

We conclude this section with a computational example illustrating the qualitative behavior of solutions. Further experiments can be found in \cite{BoNoSa18}. 
\begin{figure}[ht]
	\centering
	\includegraphics[width=0.85\textwidth]{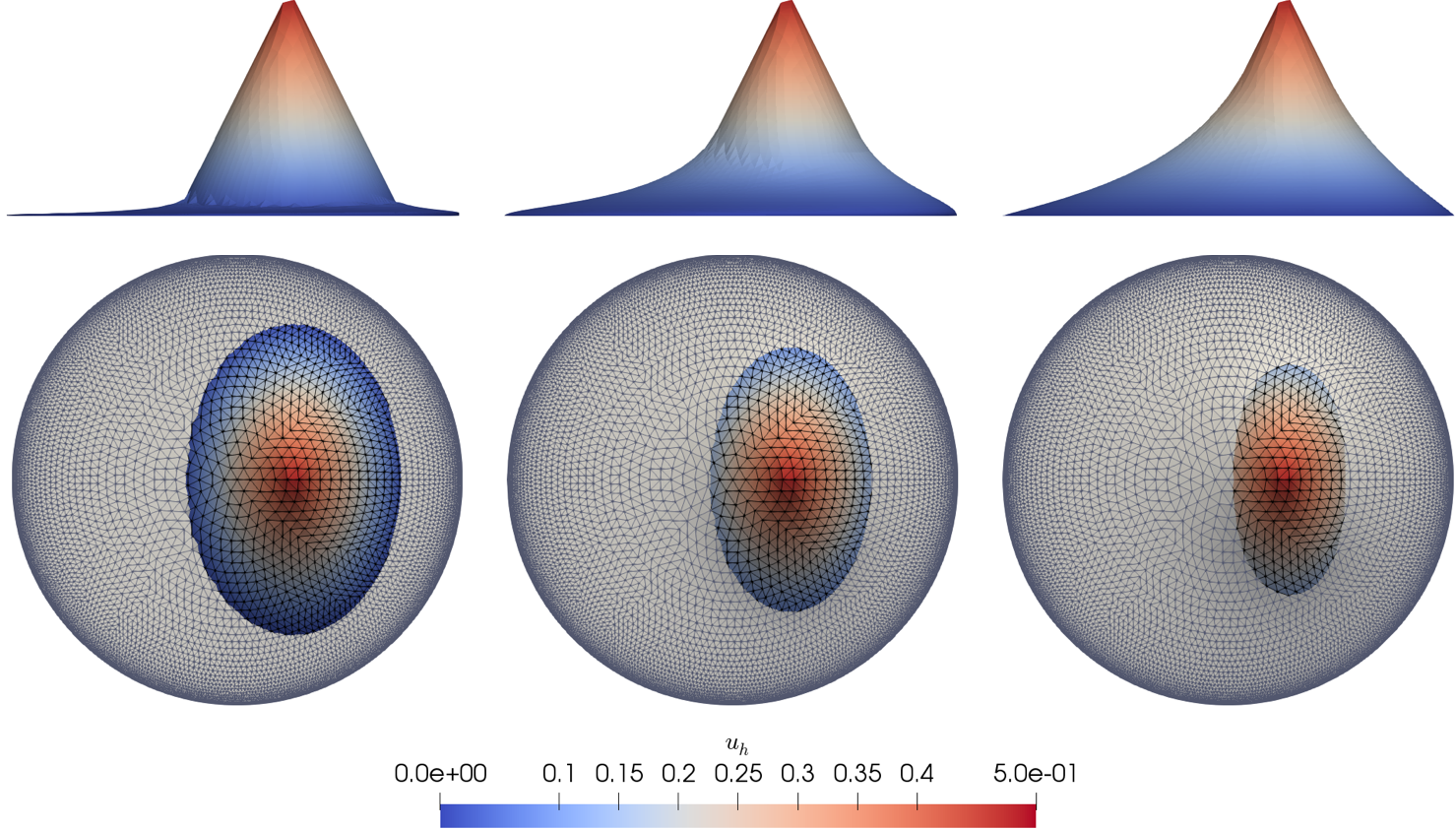}
\caption{Discrete solutions to the fractional obstacle problem for $s=0.1$ (left), $s=0.5$ (center) and $s=0.9$ (right) computed with graded meshes with $h = 2^{-5}$. Top: lateral view. Bottom: top view, with the discrete contact set highlighted.
}	\label{fig:qualitative}
\end{figure}
\begin{example}[qualitative behavior]
Consider problem \eqref{eq:obstacle}, posed in the unit ball $B_1 \subset \mathbb{R}^2$, with $f = 0$ and the obstacle 
\[ 
\chi (x_1, x_2) = \frac12 - \sqrt{ \l(x_1 - \frac14\r)^2 + \frac12 x_2^2}.
\]
\Cref{fig:qualitative} shows solutions for $s \in  \{0.1, 0.5, 0.9 \}$ on meshes graded according to \eqref{eq:H} with $\mu =2$. The coincidence set $\Lambda$, which contains a neighborhood of the singular point $(1/4,0)$ is displayed in color in the bottom view. It can be observed that, while for $s=0.9$ the discrete solution resembles what is expected for the classical obstacle problem, the solution for $s=0.1$ is much flatter in the non-coincidence set $N$. Because $f \ge 0$, the solution $u$ satisfies $u \ge 0$. Therefore, the solution $u$ approaches $\chi_+$ in the limit $s \to 0$, while $u$ is expected to touch the obstacle only at the singular point and detach immediately in the diffusion limit $s=1$. The discrete nonlinear system has been solved using a semi-smooth Newton method.
\end{example}
%

%% file: minimal_surfaces.tex
\newcommand{\wl}[1]{{\color{blue}#1}}

\section{Fractional minimal graphs} \label{sec:NMS}
In this section we discuss the fractional minimal graph problem. The line of study of this nonlinear fractional problem, that can be regarded as a nonlocal version of the classical Plateau problem, and other related problems, began with the seminal works by Imbert \cite{Imbert09} and Caffarelli, Roquejoffre and Savin \cite{CaRoSa10}. 

As a motivation for the notion of fractional minimal sets, we show how the fractional perimeter arises in the study of a nonlocal version of the Ginzburg-Landau energy, extending a well-known result for classical minimal sets \cite{ModicaMortola77}. Let $\Omega \subset \mRd$ be a bounded set with Lipschitz boundary, $\varepsilon > 0$ and define the energy
\[
\mathcal{J}_{\ve}[u;\Omega] := \frac{\ve}{2} \int_\Omega |\nabla u(x)|^2 \; dx + \frac{1}{\ve}\int_{\Omega} W(u(x)) \; dx,
\]
where $W(t) = \frac14(1-t^2)^2$ is a double-well potential. Then, for every sequence $\{ u_\eps \}$ of minimizers of $\mathcal{J}_{\ve}[u;\Omega]$ with uniformly bounded energies there exists a subsequence $\{ u_{\eps_k} \}$ such that
\[ u_{\eps_k} \to \chi_E - \chi_{E^c} \quad \mbox{in } L^1(\Omega),
\]
where $E$ is a set with minimal perimeter in $\Omega$.
We now consider a different regularization term: given $s \in (0,1/2)$, we set
\[
\mathcal{J}^s_{\ve}[u;\Omega] := \frac{1}{2}\iint_{Q_{\Omega}} \frac{|u(x) - u(y)|^2}{|x-y|^{n+2s}} \; dx dy + \frac{1}{\ve^{2s}} \int_{\Omega} W(u(x)) \; dx,
\]
where ${Q_{\Omega} = \l( \mRd \times \mRd \r) \setminus \l( {\Omega}^c \times {\Omega}^c \r)}$ as in \eqref{eq:defofinnerprod}.
The first term in the definition of $\mathcal{J}^s_\eps$ involves the $H^s(\mRd)$-norm of $u$, except that the interactions over $\Omega^c \times \Omega^c$ are removed; for a minimization problem in $\Omega$, these are indeed fixed. As proved in \cite{SaVa12Gamma}, for every sequence $\{ u_\eps \}$ of minimizers of $\mathcal{J}^s_{\ve}$ with uniformly bounded energies there exists a subsequence $\{ u_{\eps_k} \}$ such that
\[ u_{\eps_k} \to \chi_E - \chi_{E^c} \quad \mbox{in } L^1(\Omega) \quad \mbox{as } \ve_k \to 0^+.
\]
However, instead of minimizing the perimeter in $\Omega$, here the set $E$ is a $s$-minimal set in $\Omega$, because it minimizes the so-called fractional perimeter $P_s(E,\Omega)$ among all measurable sets $F \subset \mRd$ such that $F \setminus \Omega = E \setminus \Omega$. In \cite{CaRoSa10}
this notion of fractional perimeter (also known as nonlocal perimeter) was proposed, and nonlocal minimal set problems were studied.
We refer to \cite[Chapter 6]{BucurValdinoci16} and \cite{CoFi17} for nice introductory expositions to the topic and applications.


\subsection{Formulation of the problem and regularity}

Our goal is to compute fractional minimal graphs, that is, to study the nonlocal minimal surface problem under the restriction of the domain being a cylinder. Concretely, from now on we consider $\Omega' = \Omega \times \mR$ with $\Omega \subset \mRd$ being a bounded Lipschitz domain. We assume that the exterior datum is the subgraph of some uniformly bounded function $g: \mRd \setminus \Omega \to \mR$,
\begin{equation} \label{E:Def-E0}
E_0 := \l\{ (x', x_{d+1}) \colon x_{d+1} < g(x'), \; x' \in \mRd \setminus \Omega \r\}.
\end{equation}
The fractional minimal graph problem consists in finding a {\em locally} $s$-minimal set $E$ in $\Omega'$ such that $E \setminus \Omega' = E_0$. We refer to  \cite{Lomb16Approx} for details on why the notion of locally $s$-minimality is the `correct' one. Under the conditions described above, it can be shown that minimal sets need to be subgraphs, that is,
\begin{equation}\label{E:Def-E}
E \cap \Omega' = \l\{ (x', x_{d+1}) \colon x_{d+1} < u(x'), \; x' \in \Omega \r\}
\end{equation}
for some function $u$ (cf. \cite[Theorem 4.1.10]{Lombardini-thesis}).
We shall refer to such a set $E$ as a {\em nonlocal minimal graph} in $\Omega$. 

In order to find nonlocal minimal graphs, we introduce the space 
\[
\mathbb{V}^g := \{ v \colon \mRd \to \mR \; \colon \; v\big|_\Omega \in W^{2s}_1(\Omega), \ v = g \text{ in } {\Omega}^c\}
\]
(we write $\mathbb{V}^0$ whenever $g \equiv 0$)  and, considering the weight function $F_s \colon \mR \to \mR$,
\begin{equation} \label{E:def_Fs}
F_s(\rho) := \int_0^\rho \frac{\rho-r}{\l( 1+r^2\r)^{(d+1+2s)/2}} dr.
\end{equation}
we define the energy functional
\begin{equation}\label{E:NMS-Energy-Graph}
I_s[u] := \iint_{Q_{\Omega}} F_s\l(\frac{u(x)-u(y)}{|x-y|}\r) \frac{1}{|x-y|^{d+2s-1}} \;dxdy.
\end{equation}
In \cite[Chapter 4]{Lombardini-thesis} it is shown that finding nonlocal minimal graphs is equivalent to minimizing the energy $I_s$ over $\mathbb{V}^g$.
Existence of solution $u$ follows from the existence of locally $s$-minimal sets \cite{Lomb16Approx}, while uniqueness is a consequence of $I_s$ being strictly convex. We also point out that for any function $v \colon \mRd \to \mR$, its energy $I_s[v]$ is closely related to certain ${W^{2s}_1}$-seminorms \cite[Lemma 2.5]{BortLiNoch2019NMG-convergence}:
\begin{equation} \label{eq:norm_bound} \begin{aligned}
& |v|_{W^{2s}_1(\Omega)} \le C_1 + C_2 I_s[v],
& I_s[v] \le C_3 \iint_{Q_{\Omega}}  \frac{|v(x)-v(y)|}{|x-y|^{d+2s}} dxdy .
\end{aligned} \end{equation}

To give a clearer picture of the nonlocal minimal graph problem, we compare it to its classical counterpart. Given a bounded domain $\Omega \subset \mRd$ with sufficiently smooth boundary, and $g \colon \partial \Omega \to \mRd$, the classical Plateau problem consists in finding $u \colon \Omega \to \mRd$ that minimizes the graph surface area functional
\begin{equation} \label{E:MS-Energy-Graph}
I [u] := \int_\Omega \sqrt{1 + |\nabla u (x)|^2 } \, dx 
\end{equation}
among those functions $u \in H^1(\Omega)$ satisfying $u = g$ on $\partial \Omega$. By taking first variation of $I$, it follows that the minimizer $u$ satisfies
\begin{equation} \label{eq:classical-MS}
\int_\Omega \frac{\nabla u(x) \cdot \nabla v (x)}{\sqrt{1 + |\nabla u (x)|^2 }} \, dx = 0  \quad \forall \, v \in H^1_0(\Omega).
\end{equation}
The left hand side in \eqref{eq:classical-MS} consists of an $H^1$-inner product between $u$ and $v$, with a possibly degenerate weight that depends on $u$. For the nonlocal problem, after taking first variation of $I_s$ in \eqref{E:NMS-Energy-Graph}, we obtain that $u$ is a minimizer if and only if
\begin{equation}\label{E:WeakForm-NMS-Graph}
a_u(u,v) = 0 \quad \mbox{ for all }  v \in \mathbb{V}^0,
\end{equation}
where the bilinear form $a_u \colon \mathbb{V}^g \times \mathbb{V}^0 \to \mathbb{R}$ is given by
\begin{equation} \label{E:def-a}
a_u(w,v) := \iint_{Q_{\Omega}} \Gts\l(\frac{u(x)-u(y)}{|x-y|}\r) \frac{(w(x)-w(y))(v(x)-v(y))}{|x-y|^{d+1+2s}}dx dy, 
\end{equation}
and $\Gts(\rho) = \int_0^1 (1+ \rho^2 r^2)^{-(d+1+2s)/2} dr$ and hence it satisfies $\rho \Gts(\rho) = G_s(\rho) = F'_s(\rho)$. Similar to \eqref{eq:classical-MS}, the left hand side $a_u(u,v)$ in \eqref{E:WeakForm-NMS-Graph} is a weighted $H^{s+\frac{1}{2}}$-inner product with a possibly degenerate weight depending on $u$.

As for the regularity of nonlocal minimal graphs, the following result is stated in \cite[Theorem 1.1]{CabreCozzi2017gradient} and builds on the arguments in \cite{Barrios2014bootstrap, Figalli2017regularity}.

\begin{theorem}[interior smoothness of nonlocal minimal graphs] \label{thm:smoothness}
Assume $E \subset \mathbb{R}^{d+1}$ is a locally $s$-minimal set in $\Omega' = \Omega \times \mathbb{R}$, given by the subgraph of a measurable function $u$ that is bounded in an open set $\Lambda \supset \Omega$. Then, $u \in C^\infty (\Omega)$.
\end{theorem}

\begin{remark}[stickiness]\label{rem:stickiness}
\Cref{thm:smoothness} does not address boundary regularity. By using \eqref{eq:norm_bound}, it can be easily proved that $u \in W^{2s}_1(\Omega)$ but, because $2s < 1$, this does not even guarantee that $u$ has a trace on $\pp\Omega$. In fact, nonlocal minimal graphs can develop discontinuities across $\partial\Omega$. Furthermore, nonlocal minimal graphs generically exhibit this {\em sticky} behavior \cite{DipiSavinVald17, DipiSavinVald19}.
\end{remark}

\subsection{Finite element discretization}
In this section, we review the finite element discretization of the nonlocal minimal graph problem proposed in \cite{BortLiNoch2019NMG-convergence} and discuss its convergence and error estimates. 

For simplicity, we assume that $\mbox{supp}(g) \subset \Lambda$ for some bounded set $\Lambda$. As before, we take a family $\{\Th \}_{h>0}$ of conforming, simplicial and shape-regular meshes on $\Lambda$, which we impose to mesh $\Omega$ exactly. To account for non-zero boundary data, we make a slight modification on \eqref{eq:FE_space} to define the discrete spaces
\[
\mathbb{V}_h :=  \{ v \in C(\Lambda) \colon v|_T \in \mathcal{P}_1 \; \forall T \in \Th \},
\]
and we define
\[
\mathbb{V}_h^g := \{ v_h \in \mathbb{V}_h \colon \ v_h|_{\Lambda \setminus \Omega} = \Cl^c g\}, \quad
\mathbb{V}_h^0 := \{ v_h \in \mathbb{V}_h \colon \ v_h|_{\Lambda \setminus \Omega} = 0\},
\]
where $\Cl^c$ denotes the Cl\'ement interpolation operator in $\Omega^c$.

With the notation introduced above, the discrete problem seeks  $u_h \in \mathbb{V}^g_h$ such that
\begin{equation}\label{E:WeakForm-discrete}
a_{u_h}(u_h, v_h) = 0 \quad \mbox{for all } v_h \in \mathbb{V}^0_h. 
\end{equation}
Existence and uniqueness of solutions to this discrete problem follow directly from \eqref{eq:norm_bound} and the strict convexity of $I_s$. To prove the convergence of the finite element scheme, the approach in \cite{BortLiNoch2019NMG-convergence} consists in proving that the discrete energy is consistent and afterwards using a compactness argument.

\begin{theorem}[convergence for the nonlocal minimal graph problem]  \label{thm:consistency}
Let $s \in (0,1/2)$, $\Omega$ be a bounded Lipschitz domain and $g$ be uniformly bounded and satisfying $\mbox{supp}(g) \subset \Lambda$ for some bounded set $\Lambda$. Let $u$ and $u_h$ be, respectively, the solutions to \eqref{E:WeakForm-NMS-Graph} and \eqref{E:WeakForm-discrete}. Then, it holds that
\[
    \lim_{h \to 0} I_s[u_h] = I_s[u] \quad
\mbox{ and } \quad
\lim_{h \to 0} \| u - u_h \|_{W^{2r}_1(\Omega)} = 0 \quad \forall r \in [0,s).
\]
\end{theorem}

The theorem above has the important feature of guaranteeing convergence without any regularity assumption on the solution. However, it does not offer convergence rates.  We now show estimates for a geometric notion of error that mimics the one analyzed in \cite{FierroVeeser03} for the classical Plateau problem (see also \cite{BaMoNo04, DeDzEl05}). Such a notion of error is given by
\begin{equation}\label{eq:def-e}
\begin{aligned}
e^2(u,u_h) & :=\int_{\Omega} \ \Big| \widehat{\nu}(\nabla u) - \widehat{\nu}(\nabla u_h) \Big|^2 \;\frac{Q(\nabla u) + Q(\nabla u_h)}{2} \ dx , \\
& = \int_{\Omega} \ \Big( \widehat{\nu}(\nabla u) - \widehat{\nu}(\nabla u_h) \Big) \cdot \ \Big( \nabla (u-u_h), 0 \Big) dx ,
\end{aligned}
\end{equation}
where $Q(\pmb{a}) = \sqrt{1+|\pmb{a}|^2}$, $\widehat{\nu}(\pmb{a}) = \frac{(\pmb{a},-1)}{Q(\pmb{a})}$.
Because $\widehat{\nu}(\nabla u)$ is the normal unit vector on the graph of $u$, the quantity $e(u,u_h)$ is a weighted $L^2$-discrepancy between the normal vectors. 
For the nonlocal minimal graph problem, \cite{BortLiNoch2019NMG-convergence} introduced
\begin{equation}\label{eq:def-es} 
{
\begin{aligned}
e_s(u,u_h) := \l( \widetilde C_{d,s} \iint_{Q_{\Omega}} \Big( G_s\l(d_u(x,y)\r) - G_s\l(d_{u_h}(x,y)\r) \Big) \frac{d_{u-u_h}(x,y)}{|x-y|^{d-1+2s}} dxdy \r)^{1/2} ,
\end{aligned}}
\end{equation}
where $G_s(\rho) = F_s'(\rho)$, the constant $\widetilde C_{d,s} = \frac{1 - 2s}{\alpha_{d}}$, $\alpha_{d}$ is the volume of the $d$-dimensional unit ball and $d_v$ is the difference quotient of the function $v$, 
\begin{equation}\label{eq:def-d_u}
d_v(x,y) := \frac{v(x)-v(y)}{|x-y|}.
\end{equation}
In \cite{BortLiNoch2019NMG-convergence}, this novel quantity $e_s(u,u_h)$ is shown to be connected with a notion of nonlocal normal vector, and its asymptotic behavior as $s\to 1/2^-$ is established.

\begin{theorem}[asymptotics of $e_s$] \label{Thm:asymptotics-es}
For all $u,v \in H^1_0(\Lambda)$, we have
\[
\lim_{s \to {\frac{1}{2}}^-} e_s(u,v) = e(u,v).
\]
\end{theorem}

A simple `Galerkin orthogonality' argument allows to derive an error estimate for $e_s(u,u_h)$ (cf. \cite[Theorem 5.1]{BortLiNoch2019NMG-convergence}).

\begin{theorem}[geometric error]\label{thm:geometric-error}
Under the same hypothesis as in \Cref{thm:consistency}, it holds that
\begin{equation} \label{eq:geometric-error} \begin{aligned}
e_s(u,u_h) &\le C (d,s) \, \inf_{v_h \in \mathbb{V}_h^g}  \l( \iint_{Q_{\Omega}} \frac{|(u-v_h)(x)-(u-v_h)(y)|}{|x-y|^{d+2s}} dxdy \r) ^{1/2}.
\end{aligned} \end{equation}
\end{theorem}

Therefore, to obtain convergence rates with respect to $e_s(u,u_h)$, it suffices to prove interpolation estimates for the nonlocal minimizer. Although minimal graphs are expected to be discontinuous across the boundary, we still expect that $u \in BV(\Lambda)$ in general. Under this circumstance, the error estimate \eqref{eq:geometric-error}  leads to
\[
e_s(u,u_h) \le C(d,s) \, h^{1/2-s} |u|^{1/2}_{BV(\Lambda)}.
\]

\subsection{Numerical experiments}
We conclude by presenting several numerical experiments and discussing the behavior of nonlocal minimal graphs. The first example we compute is on a one-dimensional domain, and is proposed and theoretically studied in \cite[Theorem 1.2]{DipiSavinVald17} as an illustration of stickiness phenomena.

\begin{example}[stickiness in $1$D]\label{Ex:1d_stickiness}
Let $\Omega = (-1,1)$ and $g(x) = \textrm{sign}(x)$ for $x \in \Omega^c$. Discrete nonlocal minimal graphs for $s\in\{0.1, 0.25,0.4\}$ are shown in \Cref{F:Ex-stickiness} (left).

\begin{figure}[!htb]
	\begin{center}
		\includegraphics[width=0.4\linewidth]{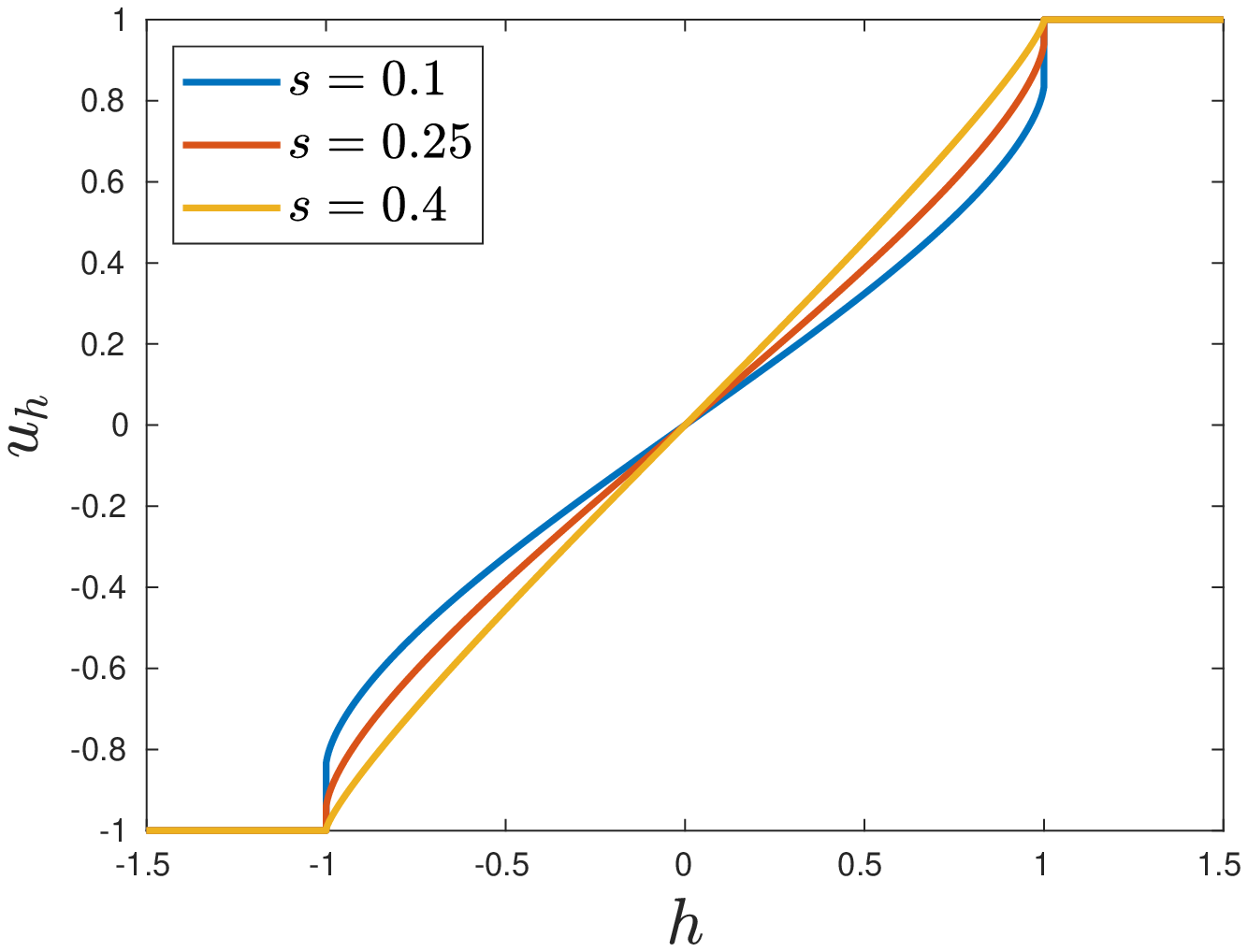}
		\includegraphics[width=0.4\linewidth]{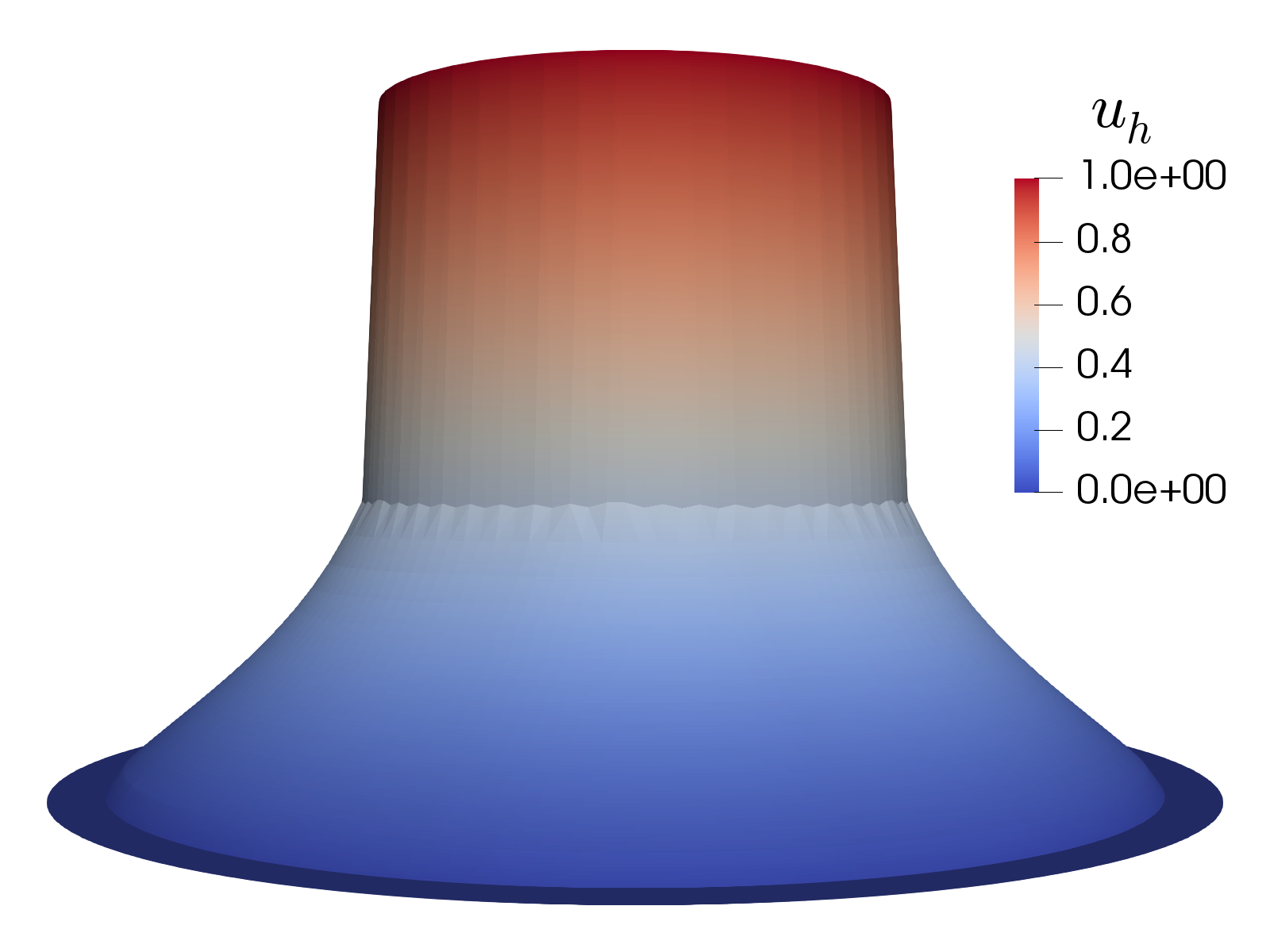}
	\end{center}
	\caption{\small Left: plot of $u_h$ for $s = 0.1, 0.25, 0.4$ (from left to right) in \Cref{Ex:1d_stickiness}. Right: plot of $u_h$ for uniform meshes with $h = 2^{-5}$ in \Cref{Ex:2d_annulus}. }
	\label{F:Ex-stickiness}
\end{figure}
Although our method requires discrete functions to be continuous across the boundary of $\Omega$, the presented $1$D picture clearly suggests a stickiness phenomenon.
In addition, the plot also indicates that stickiness becomes more observable when $s$ gets closer to $0$.  
\end{example}

\begin{example}[stickiness in an annulus]\label{Ex:2d_annulus}
Let $\Omega = B_1 \setminus \overline{B}_{1/2} \subset \mR^2$, where $B_r$ denotes an open ball with radius $r$ centered at the origin, and let $g = \chi_{B_{1/2}}$.
The discrete nonlocal minimal graph for $s = 0.25$ is plotted in \Cref{F:Ex-stickiness} (right).
In this example, stickiness is clearly observed on $\partial B_{1/2}$, while the stickiness on $\partial B_1$ is less noticeable.
\end{example}


\begin{example}[effect of $s$]\label{Ex:2d_circle}
Let $\Omega = B_1 \subset \mR^2$, $g = \chi_{B_{3/2}\setminus\Omega}$. \Cref{F:NMS-Ex_multis} shows minimizers for several values of $s$.
As $s \to 1/2$, the nonlocal minimal graphs get closer to the classical minimal graph, which is trivially constant in $\Omega$. On the other hand, as $s$ decreases we observe a stronger jump across $\pp\Omega$.

\begin{figure}[!htb]
	\begin{center}
		\includegraphics[width=0.85\linewidth]{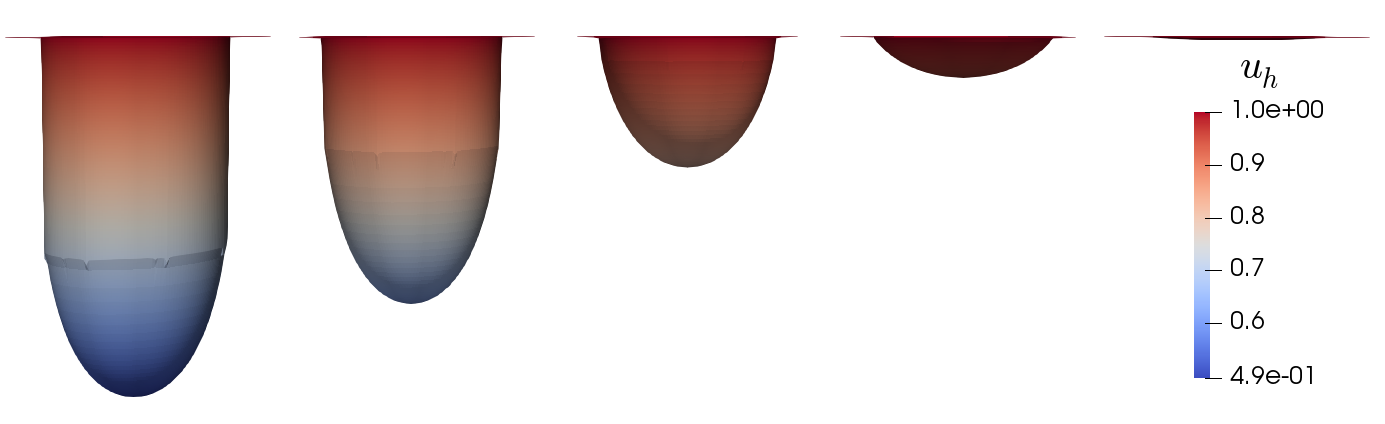}
	\end{center}
	\caption{\small Plot of $u_h$ for $s = 0.01, 0.1, 0.25, 0.4, 0.49$ (from left to right) with uniform $h = 2^{-5}$ in \Cref{Ex:2d_circle}.}
	\label{F:NMS-Ex_multis}
\end{figure}
\end{example}

\wl{In the end, we make brief comments on the computational side for the examples presented above. In all of the experiments for fractional minimal graphs, we use Newton's method to solve the nonlinear equation \eqref{E:WeakForm-discrete}. Although the Jacobian matrix $\vA_{u_k}$ in the iterative process corresponds to an $H^{s+1/2}-$inner product $a_{u_k} (w,v)$ with degenerate weights, our experiments indicate its condition number behaves like 
	$$
	\kappa(\vA_{u_k}) \approx \mathcal{O}\left( N^{2(s+\frac{1}{2})/d} \right) 
	$$ 
for $u_k$ whose gradient blows up near the boundary $\pO$, quasi-uniform meshes and dimensions $d=1,2$. This behavior is the same as the one in linear fractional diffusion with order $s+\frac{1}{2}$. However, the degenerate weight does bring more difficulties in preconditioning.

Another thing worth to be pointed out is the use of Dirichlet condition for our discrete space $\mathbb{V}_h^g$ requires the discrete function $u_h$ to be continuous across $\partial \Omega$. Due to the stickiness phenomenon in \Cref{rem:stickiness}, this may not be true for the solution $u$ of the minimal graph problem. Fortunately, this does not preclude the convergence in `trace blind' fractional Sobolev spaces $W^{2s}_1(\Omega)$, and we are still able to capture discontinuities across the boundary in practice. While permitting discontinuities would be desirable, it has conflicts with using Newton's method in solving \eqref{E:WeakForm-discrete} because the bilinear form $a_u(w,v)$ in \eqref{E:def-a} may not be well-defined. The question of how to solve the nonlinear equation \eqref{E:WeakForm-discrete} faster when allowing discontinuous across $\partial \Omega$ is still under investigation.
} 
